\documentclass[12pt ]{article}

\newcommand{\bc}{\begin{center}}
\newcommand{\ec}{\end{center}}
\newcommand{\ba}{\begin{array}}
\newcommand{\ea}{\end{array}}

\date{}

\newtheorem{theorem}{Theorem}[section]
\newtheorem{corollary}{Corollary}[section]

\newcommand{\na}{\nabla}

\newcommand{\vob}{V(x_{0},\,r)}
\newcommand{\vk}{V(x_{0},\,kr)}
\newcommand{\vr}{V_{0}(r)}
\newcommand{\cm}{C_{0}^{\infty}(M)}
\newcommand{\rn}{R^{n}}

\newcommand{\py}{\frac{1}{p}}
\newcommand{\qy}{\frac{1}{q}}
\newcommand{\ny}{\frac{1}{n}}

\title{ General Sobolev Inequality on  Riemannian Manifold  \ $^{\ast}$
 \footnotetext{$\ast$   Project partially supported by NFSC No. 10271089;\,\, and partially supported by Fund of Department of education of Fujian Province
; No.JA04266.  } }
\author{ Qihua Ruan$^{1}$  \and Zhihua
Chen$^{2}$}
\begin{document}

\maketitle
\begin {center}
{ 1. Department of Applied Mathematics,Tongji University,
Shanghai,\,200092 P.R.China;

 Department of Mathematics, Putian College, Putian, Fujian  351100
P.R.China.\,\,\,\,\,\,

Email:\,\,ruanqihua@163.com

 2.  Department of Applied
Mathematics,Tongji University, Shanghai,\,200092 P.R.China.

Email:\,\,zzzhhc@tongji.edu.cn }
\end{center}
\thispagestyle{empty}

\begin{abstract}  
Let M be a complete n-dimensional Riemannian manifold, if the
sobolev inqualities hold on M, then the geodesic ball has maximal
volume growth; if  the Ricci curvature of M is nonnegative, and
one of the  general Sobolev inequalities holds on M, then M is
diffeomorphic to $R^{n}$.

\vspace{0.3cm} \noindent {\bf Keywords:}   General Sobolev
inequality, Ricci curvature, diffeomorphic.

\vspace{0.3cm} \noindent {\bf Mathematics Subject Classification
(2000):}\ 53C20, 53C21.
\end{abstract}

\setcounter{equation}{0}
\section{Introduction}
Let M be any complete n-dimensional ($n \geq 2$ ) Riemannian
manifold with nonnegative Ricci curvature,  $\cm$ be the space of
smooth functions with compact support in M. Denote by dv and $\na$
the Riemannian volume element and the gradient operator of M,
respectively.

It is well known in [1] that Ledoux showed that:  If one of the
following Sobolev inequalities is satisfied,
$$||f||_{p}\leq
C_{0}||\na f||_{q},\forall f \in \cm,\,\, 1 \leq q <n,\,\,\py =
\qy - \ny.\leqno(1)$$ where $C_{0} $ is the optimal constant in
$\rn,$  denote $||f||_{p}$ by the $ L^{p}$ norm of function $f$;
then M is \ isometric to $\rn.$

The basic idea of Ledoux's result is to find a function in $\cm$,
 then one can substitute it to (1) and obtain that $\vob \geq\vr$,
here $\vob$ denote the volume of the geodesic ball $B(x_{0},r)$ of
radius $r$ with center $x_{0}$, and $\vr$ the volume of the
Euclidean ball of radius $r$ in $\rn.$  Since the Ricci curvature
of M is nonnagtive , from Bishop's comparison theorem[2], we know
that $\vob \leq \vr,$  so M is isometric to $\rn.$

Later Xia combined Ledoux's method with Cheeger and Colding's
result[3], which is that given an integer $n\geq 2$, there exists
a constant $\delta (n)>0$ such that any n-dimensional complete
Riemannian manifold with nonnegative Ricci curvature and $\vob
\geq(1-\delta(n))\vr$ for some $x_{0}\in M$ and all $ r>0 $ is
diffeomorphic to $\rn$; and proved the following theorem [4]:
\begin{theorem}Let M be a complete n-dimensional Riemannian manifold with
nonnegative Ricci curvature, if one of the following sobolev
inequalities is satisfied,
$$||f||_{p}\leq
C_{1}||\na f||_{q},\forall f \in \cm,\,\, 1 \leq q <n,\,\,\py =
\qy - \ny.\leqno(2)$$ where the positive constant $C_{1}>C_{0}$,
then M is diffiomorphic to $\rn$.
\end{theorem}
As we known, the Sobolev inequality(2) belongs to a general
Sobolev  inequalities of the type
$$||f||_{t}\leq
C||f||_{s}^{\theta}\ ||\na f||_{q}^{1-\theta},\forall f \in
\cm,\,\,
\frac{1}{t}=\frac{\theta}{s}+\frac{1-\theta}{p}.\leqno(3)$$ (see
[5]), where C is a positive constant. Inequality (2) corresponds
to $\theta = 0$. When $q=r=2$, and $\theta=2/(n+2)$, it
corresponds to the Nash inequality
$$\left(\int\mid f \mid^{2}dv\right)^{1+\frac{2}{n}} \leq C\left(\int\mid f \mid dv\right)
^{\frac{4}{n}} \int\mid \nabla f \mid^{2}dv,\,\forall f \in
C_{0}^{\infty}(M),\leqno(4)$$ [ see(6)], where C is a positive
constant.

In \cite{}, Ledoux conjecture two problem: one is that if one of
Nash inequalities is satisfied on M, then M is isometric to $\rn$;
the other is that without any curvature assumption, if the Sobolev
inequalities are satisfied on M, then the geodesic ball has
maximal volume growth, namely $\lim\limits_{r\rightarrow +\infty}
\frac{\vob}{\vr}\geq C$, where positive constant C depends only on
$n$ and $q$. The first problem was proved by Druet, Hebey, and
Vaugon(see \cite{}). In this paper, we want to confirm the other
problem, that is
\begin{theorem}Let M be a complete n-dimensional Riemannian
manifold, if the Sobolev inequalities (2) are satisfied on M, then
the geodesic ball in M has maximal volume growth.
\end{theorem}
We also generalize Theorem 1.1 to  obtain the following result:
\begin{theorem}
Let M be a complete n-dimensional Riemannian manifold with
nonnegative Ricci curvature, if one of the general Sobolev
inequalities (3) is satisfied on M, then M is diffeomorphic to
$\rn$.
\end{theorem}

\setcounter{equation}{0}
\section{Proof of  Theorem 1.2 and 1.3 }
Now we first prove Theorem 1.2.
\\ {\bf \large Proof:} As we known,
there always exists a cut-off function $ f_{m}\in \cm$  such that
$f_{m}\equiv 1$ on $B(x_{0},2^{m}r)$, $f_{m}\equiv 0$ outside
$B(x_{0},2^{m+1}r)$, $|\nabla f_{m}|\leq \frac{1}{2^{m}r}$ on
$B(x_{0},2^{m+1}r)\setminus B(x_{0},2^{m}r)$ for any real integer
$m$. Substituting them to (1), we have that
$$V(x_{0}, 2^{m}r)^{\frac{1}{p}}\leq \frac{C_{0}}{2^{m}r}V(x_{0}, 2^{m+1}r)^{\frac{1}{q}}$$
Thus set $\frac{q}{p}=a $, we get that
$$V(x_{0}, r)\geq
(\frac{r}{2C_{0}})^{q}V(x_{0}, \frac{r}{2})^{a}\hspace{6cm}$$
$$\geq (C_{0}^{-1}r)^{q(1+a)}2^{-q(1+2a)}V(x_{0}, \frac{r}{2^{2}})^{a^{2}}\hspace{1.9cm}$$
$$\geq (C_{0}^{-1}r)^{q(1+a+a^{2})}2^{-q(1+2a+3a^{2})}V(x_{0}, \frac{r}{2^{3}})
^{a^{3}}$$
$$\hspace{0.1cm}\geq\cdot\cdot\cdot\geq (C_{0}^{-1}r)^{q\sum\limits_{i=0}^{k}a^{i}}2^{-q\sum\limits_{i=0}
^{k}(i+1)a^{i}}V(x_{0}, \frac{r}{2^{k+1}}) ^{a^{k+1}}$$
$$\geq\cdot\cdot\cdot\geq (C_{0}^{-1}r)^{q\sum\limits_{i=0}^{k}a^{i}}2^{-q\sum\limits_{i=0}
^{k}(i+1)a^{i}}V(x_{0}, 1) ^{a^{k+1}} \ (r\geq 2^{k+1})$$ Let
$k\rightarrow +\infty,\ r\rightarrow +\infty$, since
$a=\frac{q}{p}<1$, then we see that
$$\lim\limits_{r\rightarrow +\infty}\frac{V(x_{0}, r)}{r^{n}}\geq C_{0}^{-n}2^{n-\frac{n^{2}}{q}}$$
So we proved Theorem 1.2.
\\ Next we prove Theorem 1.3
\\ {\bf \large Proof:} As proof of Theorem 1.2, there always exists
a cut-off function $f\in \cm$  such that $f\equiv 1$ on
$B(x_{0},r)$, $f\equiv 0$ outside $B(x_{0},kr)$, $|\nabla f|\leq
\frac{1}{(k-1)r}$ on $B(x_{0},kr)\setminus B(x_{0},r)$ for some
$k>1$.

Substituting it to (1), we have that
$$\vob^{\frac{1}{t}}\leq C \vk^{\frac{\theta}{s}}(\frac{1}{(k-1)r})^{1-\theta}
\vk^{\frac{1-\theta}{q}}$$
$$\vob^{\frac{1}{t}}\leq \frac{C}{((k-1)r)^{1-\theta}}\vk^{\frac{\theta}{s}+
\frac{1-\theta}{q}}\hspace{2cm}$$
$$\vob^{\frac{1}{t}}\leq \frac{C}{((k-1)r)^{1-\theta}}k^{n(\frac{\theta}{s}+
\frac{1-\theta}{q})} \vob^{\frac{\theta}{s}+
\frac{1-\theta}{q}}\hspace{0.8cm}$$
$$1\leq \frac{Ck^{n(\frac{\theta}{s}+
\frac{1-\theta}{q})}}{((k-1)r)^{1-\theta}}\vob^{\frac{\theta}{s}+
\frac{1-\theta}{q}-\frac{1}{t}}\hspace{0.7cm}$$
$$1\leq \frac{Ck^{n(\frac{\theta}{s}+
\frac{1-\theta}{q})}}{((k-1)r)^{1-\theta}}\vob^{\frac{1-\theta}{n}}\hspace{1.5cm}$$
$$\vob\geq \frac{(k-1)^{n}}{C^{\frac{n}{1-\theta}}k^{\frac{n^{2}\theta}{s(1-\theta)}+
\frac{n^{2}}{q}}}r^{n}\hspace{4.1cm}$$ For fixed number $\theta$
and $q$, we can choose one
$k>(C^{-\frac{1}{1-\theta}}\omega_{n}^{-\frac{1}{n}})^{\frac{1}{\frac{n\theta}{s(1-\theta)}+
\frac{n}{q}-1}}$ such that $\vob \geq C'\vr$ and $C' <1$. Then
there always exists a $\delta(n)> 0$ such that $\vob \geq
(1-\delta(n)) \vr$. Thus from Cheeger and Colding's result, we see
that M is diffeomorphic to $\rn$

If $\theta=0$, the general Sobolev inequality change into the
Sobolev inequality, thus we generalize Xia's result. If $q=r=2$
and $\theta=2/(n+2)$, then the general Sobolev inequality become
the Nash inequality, so we get the following corollary.
\begin{corollary}
Let M be a complete n-dimensional Riemannian manifold with
nonnegative Ricci curvature, if one of the Nash inequalities (4)
is satisfied on M, then M is diffeomorphic to $\rn$.
\end{corollary}

\end{document}